\documentclass[12pt,english]{article}
\usepackage[T1]{fontenc}
\usepackage[latin9]{inputenc}
\usepackage{geometry}
\usepackage{amsmath}
\usepackage{amssymb}
\usepackage{graphicx}
\usepackage{esint}

\makeatletter
\date{}
\numberwithin{equation}{section}

\newtheorem{theorem}{Theorem}

\newcommand{\qed}{\nobreak \ifvmode \relax \else       \ifdim\lastskip<1.5em \hskip-\lastskip       \hskip1.5em plus0em minus0.5em \fi \nobreak       \vrule height0.75em width0.5em depth0.25em\fi}

\makeatother

\usepackage{babel}
\begin{document}

\title{Discriminants of Polynomials Related to Chebyshev Polynomials: The
``Mutt and Jeff'' Syndrome}

\author{Khang Tran\\
University of Illinois at Urbana-Champaign}
\maketitle
\begin{abstract}
The discriminants of certain polynomials related to Chebyshev polynomials
factor into the product of two polynomials, one of which has coefficients
that are much larger than the other's. Remarkably, these polynomials
of dissimilar size have ``almost'' the same roots, and their discriminants
involve exactly the same prime factors.
\end{abstract}

\section{Introduction}

The discriminants of the Chebyshev $T_{n}(x)$ and $U_{n}(x)$ polynomials
are given by simple and elegant formulas:
\begin{equation}
\mbox{Disc}_{x}T_{n}(x)=2^{(n-1)^{2}}n^{n}\label{eq:disT}
\end{equation}
and 
\begin{equation}
\mbox{Disc}_{x}U_{n}(x)=2^{n^{2}}(n+1)^{n-2}.\label{eq:disU}
\end{equation}
Is there anything comparable for linear combinations or integral transforms
of Chebyshev polynomials? For a special type of linear combination,
formula \eqref{eq:disU} was generalized in \cite{ds} to 
\[
\mbox{Disc}_{x}(U_{n}(x)+kU_{n-1}(x))=2^{n(n-1)}a_{n-1}(k)
\]
 where 
\[
a_{n-1}(k)=(-1)^{n}\frac{(2n+1)^{n}k^{n}}{(n+1)^{2}-n^{2}k^{2}}\left(U_{n}\left(-\frac{a+1+nk^{2}}{(2n+1)k}\right)+kU_{n-1}\left(-\frac{n+1+nk^{2}}{(2n+1)k}\right)\right).
\]
For various formulas related to this type of linear combination, see
\cite{ds} and \cite{gi}. What can be said about discriminants in
$z$ of an integral transform
\[
(Sp)(z)=\frac{1}{2}\int_{-z}^{z}p'(t)(x-t^{2})dt
\]
where $p$ is some Chebyshev polynomial? Our goal here is to show
that the resulting polynomials have discriminants that factor in a
remarkable way. More precisely, when $p=U_{2n-1}(z)$ the discriminant
factors into two polynomials, one of which has coefficients that are
much larger than the coefficients of the other. Moreover, these two
polynomials have ``almost'' the same roots, and their discriminants
involve exactly the same prime factors. For example, when $n=6$,
these two polynomials are 
\begin{eqnarray*}
M(x) & = & -143+2002x-9152x^{2}+18304x^{3}-16640x^{4}+5632x^{5}\\
J(x) & = & -2606483707+826014609706x-10410224034496x^{2}\\
 &  & +40393170792832x^{3}-60482893968640x^{4}+30616119778816x^{5}.
\end{eqnarray*}
The discriminants of $M(x)$ and $J(x)$ are $2^{64}3^{4}11^{3}13^{4}$
and $2^{40}3^{4}11^{35}13^{44}$ respectively. The roots of $M(x)$
rounded to 5 digits are 
\[
\{0.13438,0.36174,0.62420,0.85150,0.98272\}
\]
whereas those of $J(x)$ are 
\[
\{0.0032902,0.13452,0.36181,0.62428,0.85163\}.
\]
We notice that in this case, the discriminants of these two polynomials
have the same (rather small) prime factors. And except for the root
0.98272 of $M(x)$ and 0.0032902 of $J(x)$, the remaining roots are
pairwise close. In fact, for $n$ large we can show, after deleting
one root from each of $M(x)$ and $J(x)$, that the remaining roots
can be paired in such a way that the distance between any two in a
given pair is at most $1/2n^{2}$. 

We call the small polynomial $M(x)$ the ``Mutt'' polynomial and
the large polynomial $J(x)$ the ``Jeff'' polynomial after two American
comic strip characters drawn by Bud Fisher. They were an inseparable
pair, one of whom (``Mutt'') was very short compared to the other
(``Jeff''). These two names are suggested by Kenneth Stolarsky.

Now, what will happen if we take the discriminant in $t$ of $U'_{2n-1}(t)(x-t^{2})$
instead of taking the integral? It is not difficult to show that this
discriminant is
\[
Cx\left(U'_{2n-1}(\sqrt{x})\right)^{4}.
\]
where $C$ is a constant depending on $n$. As we will see, one remarkable
property about the polynomial $U'_{2n-1}(\sqrt{x})$ is that its discriminant
has the same factors as the discriminants of our polynomials $M(x)$
and $J(x)$ and its roots are almost the same as those of these two
polynomials.

The main results are given by theorems 2, 3 and 4 in sections 4, 5
and 6 respectively. In section 2 we provide the reader with a convenient
summary of notations and properties of Chebyshev polynomials, discriminants,
and resultants. Section 3 of this paper defines the Mutt and the Jeff
polynomials. Section 4 and 5 of this paper analyze the discriminants
of $J(x)$and $M(x)$. Section 6 shows that the roots of these two
polynomials are pairwise close after deleting one root from each.

\section{Discriminant, Resultant and Chebyshev polynomials}

The discriminant of a polynomial $P(x)$ of degree $n$ and leading
coefficient $\gamma$ is

\begin{equation}
\mbox{Disc}_{x}P(x)=\gamma^{2n-2}\prod_{i<j}(r_{i}-r_{j})^{2}\label{eq:defdis}
\end{equation}
where $r_{1},r_{2},\ldots,r_{n}$ are the roots of $P(x)$. The resultant
of two polynomials $P(x)$ and $Q(x)$ of degrees $n,m$ and leading
coefficients $p,q$ respectively is 
\begin{eqnarray*}
\mbox{Res}_{x}(P(x),Q(x)) & = & p^{m}\prod_{P(x_{i})=0}Q(x_{i})\\
 & = & q^{n}\prod_{Q(x_{i})=0}P(x_{i}).
\end{eqnarray*}
The discriminant of $P(x)$ can be computed in terms of the resultant
between this polynomial and its derivative 
\begin{eqnarray}
\mbox{Disc}_{x}P(x) & = & (-1)^{n(n-1)/2}\frac{1}{\gamma}\mbox{Res}(P(x),P'(x))\label{eq:disres}\\
 & = & (-1)^{n(n-1)/2}\gamma^{n-2}\prod_{i\le n}P'(r_{i}).\label{eq:disdif}
\end{eqnarray}

In this paper we compute the discriminants and resultants of polynomials
related to Chebyshev polynomials. Here we review some basic properties
of the Chebyshev polynomials. The Chebyshev polynomial of the second
kind $U_{n}(x)$ has the derivative 
\begin{equation}
U'_{n}(x)=\frac{(n+1)T_{n+1}-xU_{n}}{x^{2}-1}.\label{eq:difU}
\end{equation}

The derivative of the Chebyshev polynomial of the first kind $T_{n}(x)$
is 
\begin{equation}
T'_{n}(x)=nU_{n-1}(x).\label{eq:difT}
\end{equation}
The connections between these polynomials are given by 
\begin{eqnarray}
T_{n}(x) & = & \frac{1}{2}\left(U_{n}(x)-U_{n-2}(x)\right)\label{eq:TU1}\\
 & = & U_{n}(x)-xU_{n-1}(x)\label{eq:TU2}\\
 & = & xT_{n-1}(x)-(1-x^{2})U_{n-2}(x).\label{eq:TU3}
\end{eqnarray}
Mourad Ismail \cite{ismail} applied the following theorem from Von
J. Schur \cite{schur} to compute the generalized discriminants of
the generalized orthogonal polynomials:

\begin{theorem}

Let $p_{n}(x)$ be a sequence of polynomials satisfying the recurrence
relation
\[
p_{n}(x)=(a_{n}x+b_{n})-c_{n}p_{n-2}(x)
\]
 and the initial conditions $p_{0}(x)=1$, $p_{1}(x)=a_{1}x+b_{1}$.
Assume that $a_{1}a_{n}c_{n}\ne0$ for $n>1$. Then
\[
\prod_{p_{n}(x_{i})=0}p_{n-1}(x_{i})=(-1)^{n(n-1)/2}\prod_{j=1}^{n}a_{j}^{n-2j+1}c_{j}^{j-1}
\]
for $n\ge1$. 

\end{theorem}

Ismail's idea is that if we can construct $A_{n}(x)$ and $B_{n}(x)$
so that 
\[
p_{n}'(x)=A_{n}(x)p_{n-1}(x)+B_{n}(x)p_{n}(x).
\]
 Then 
\[
\mbox{Disc}_{x}p_{n}(x)=\gamma^{n-2}\prod_{j=1}^{n}a_{j}^{n-2j+1}c_{j}^{j-1}\prod_{p_{n}(x_{i})=0}A_{n}(x_{i}).
\]
In the special case of the Chebyshev polynomial satisfying the recurrence
relation $U_{n}(x)=2xU_{n-1}(x)-U_{n-2}(x)$ and the differental equation
\begin{eqnarray*}
U'_{n}(x) & = & \frac{2xnU_{n}-(n+1)U_{n-1}}{x^{2}-1},
\end{eqnarray*}
we obtain
\begin{equation}
\prod_{U_{2n}(x)=0}U_{2n-1}(x)=1\label{eq:resU}
\end{equation}
and 
\[
\mbox{Disc}_{x}U_{n}(x)=2^{n^{2}}(n+1)^{n-2}.
\]

For the discriminants of various classes of polynomials, see \cite{aar,apostol,ds,gkz,gi}.

\section{The Mutt and Jeff polynomial pair}

In this section we will show that the discriminant of the integral
transform of Chebyshev polynomial factors into the square of the product
of the Mutt and Jeff polynomials whose formulas will be provided.
While the formula for the Mutt polynomial can be given explicitly
in terms of the Chebyshev polynomials, we can only describe the Jeff
polynomial by its roots. By integration by parts, we have 
\begin{eqnarray}
(SU_{2n-1})(z) & = & \frac{1}{2}\int_{-z}^{z}U'_{2n-1}(t)(x-t^{2})dt\nonumber \\
 & = & (x-z^{2})U_{2n-1}(z)+\int_{-z}^{z}U_{2n-1}(t)tdt\nonumber \\
 & = & (x-z^{2})U_{2n-1}(z)+\frac{1}{2}\int_{-z}^{z}U_{2n}(t)+U_{2n-2}(t)dt\nonumber \\
 & = & (x-z^{2})U_{2n-1}(z)+\left(\frac{T_{2n+1}(z)}{2n+1}+\frac{T_{2n-1}(z)}{2n-1}\right).\label{eq:IU}
\end{eqnarray}

This is a polynomial in $z$ with the leading coefficient 
\[
-2^{2n-1}+\frac{2^{2n}}{2n+1}=-\frac{(2n-1)}{(2n+1)}2^{2n-1}.
\]

The discriminant of this polynomial in $z$ is
\begin{eqnarray*}
 &  & C_{n}\mbox{Res}_{z}((SU_{2n-1})(z),U'_{2n-1}(z)(x-z^{2}))\\
 & = & C_{n}\mbox{Res}_{z}((SU_{2n-1})(z),x-z^{2})\mbox{Res}_{z}((SU_{2n})(z),U'_{2n-1}(z))\\
 & = & C_{n}\left((2n-1)T_{2n+1}(\sqrt{x})+(2n+1)T_{2n-1}(\sqrt{x})\right)^{2}\left(\prod_{\substack{U'_{2n-1}(\zeta_{i})=0\\
\zeta_{i}>0
}
}(SU_{2n-1})(\zeta_{i})\right)^{2}
\end{eqnarray*}
where $C_{n}$ is a rational number depending only on $n$ and can
be different in each occurrence. Also the factors of its numerator
and denominator can only be 2 or factors of $2n-1$, $2n+1$. From
this, we define the Mutt polynomial 
\begin{eqnarray*}
M(x) & = & \frac{(2n-1)T_{2n+1}(\sqrt{x})+(2n+1)T_{2n-1}(\sqrt{x})}{x\sqrt{x}}.
\end{eqnarray*}
Also, we can define, within a plus or minus sign, the Jeff polynomial
$J(x)$$\in\mathbb{Z}[X]$ as the polynomial of degree $n-1$ whose
coefficients are relatively prime and for which
\begin{equation}
C_{n}\left(\prod_{\substack{U'_{2n-1}(\zeta_{i})=0\\
\zeta_{i}>0
}
}(SU_{2n-1})(\zeta_{i})\right)^{2}=A_{n}J^{2}(x)\label{eq:defJ}
\end{equation}
where $A_{n}$ is a suitable rational number.

\section{The Discriminant of $J(x)$ }

In this section, we find the factors of the discriminant of $J(x)$.
In particular we will prove the theorem:

\begin{theorem}

The discriminant of $J(x)$ has the same prime factors as those of
the discriminant of $U'_{2n-1}(\sqrt{x})$. Also
\[
\mbox{Disc}_{x}U'_{2n-1}(\sqrt{x})=3(2n+1)^{n-2}(2n-1)^{n-3}n^{n-3}2^{2n^{2}-3n-1}.
\]

\end{theorem}

Since we do not have an explicit formula for $J(x)$, we will compute
the discriminant as the square of the product of the distances between
the roots. We first note that 
\[
\left(\prod_{\substack{U'_{2n-1}(\zeta_{i})=0\\
\zeta_{i}>0
}
}(SU_{2n-1})(\zeta_{i})\right)^{2}
\]
is a polynomial in $x$ whose leading coefficient is 
\begin{eqnarray*}
\pm\prod_{U'_{2n-1}(\zeta_{i})=0}U_{2n-1}(\zeta_{i}) & =\pm & \frac{\mbox{Disc}_{x}U_{2n-1}(x)}{(2n-1)^{2n-1}}\\
 & =\pm & \frac{2^{(2n-1)^{2}}(2n)^{2n-3}}{(2n-1)^{2n-1}}.
\end{eqnarray*}
Thus from \eqref{eq:defJ}, the factors of the leading coefficient
of $J(x)$ can only be 2 or factors of $(2n-1)$, $(2n+1)$, $n$. 

We now consider the roots of $J(x)$. According to the formula \eqref{eq:IU},
$J(x)$ has $n-1$ real roots given by the formula 
\[
-\frac{T_{2n+1}(\zeta_{i})}{(2n+1)U_{2n-1}(\zeta_{i})}-\frac{T_{2n-1}(\zeta_{i})}{(2n-1)U_{2n-1}(\zeta_{i})}+\zeta_{i}^{2}
\]
where $U'_{2n-1}(\zeta_{i})=0$ and $\zeta_{i}>0$. The derivative
formula \eqref{eq:rootJ} implies that $2nT_{2n}(\zeta_{i})=\zeta_{i}U_{2n-1}(\zeta_{i}).$
Thus the equation \eqref{eq:TU3} gives 
\begin{eqnarray*}
T_{2n+1}(\zeta_{i}) & = & \zeta_{i}T_{2n}(\zeta_{i})-(1-\zeta_{i}^{2})U_{2n-1}(\zeta_{i})\\
 & = & U_{2n-1}(\zeta_{i})\left(\zeta_{i}^{2}(1+\frac{1}{2n})-1\right).
\end{eqnarray*}
Hence the roots of $J(x)$ can be written as below: 
\begin{eqnarray}
 &  & -\frac{T_{2n+1}(\zeta_{i})}{(2n+1)U_{2n-1}(\zeta_{i})}-\frac{T_{2n-1}(\zeta_{i})}{(2n-1)U_{2n-1}(\zeta_{i})}+\zeta_{i}^{2}\nonumber \\
 & = & \left(\frac{1}{2n-1}-\frac{1}{2n+1}\right)\frac{T_{2n+1}(\zeta_{i})}{U_{2n-1}(\zeta_{i})}-\frac{2\zeta_{i}T_{2n}(\zeta_{i})}{(2n-1)U_{2n-1}(\zeta_{i})}+\zeta_{i}^{2}\nonumber \\
 & = & \left(\frac{1}{2n-1}-\frac{1}{2n+1}\right)\left(\zeta_{i}^{2}\left(1+\frac{1}{2n}\right)-1\right)-\frac{\zeta_{i}^{2}}{n(2n-1)}+\zeta_{i}^{2}\nonumber \\
 & = & \zeta_{i}^{2}-\frac{2}{(2n+1)(2n-1)}.\label{eq:rootJ}
\end{eqnarray}
Thus from the definition of discriminant as the product of differences
between roots, it suffices to consider the discriminant of $U'_{2n-1}(\sqrt{x})$.
The formulas \eqref{eq:difU} and \eqref{eq:difT} give

\begin{eqnarray}
2\sqrt{x}U''_{2n-1}(\sqrt{x}) & = & -\frac{2\sqrt{x}}{x-1}U'_{2n-1}(\sqrt{x})+\frac{1}{x-1}\left(4n^{2}U_{2n-1}(\sqrt{x})-\sqrt{x}U'_{2n-1}(\sqrt{x})-U_{2n-1}(\sqrt{x})\right)\nonumber \\
 & = & \frac{-3\sqrt{x}}{x-1}U'_{2n-1}(\sqrt{x})+\frac{4n^{2}-1}{x-1}U_{2n-1}(\sqrt{x}).\label{eq:difj}
\end{eqnarray}
We note that $U'_{2n-1}(\sqrt{x})$ is a polynomial of degree $n-1$
with the leading coefficient 
\[
(2n-1)2^{2n-1}.
\]

From the definition of discriminant in terms of resultant \eqref{eq:disdif}
and the formulas \eqref{eq:disU} and \eqref{eq:difj}, we have
\begin{eqnarray*}
\mbox{Disc}_{x}U'_{2n-1}(\sqrt{x}) & = & (2n-1)^{n-3}2^{(2n-1)(n-3)}\prod_{U'_{2n-1}(\sqrt{x_{i}})=0}\frac{4n^{2}-1}{x_{i}-1}\frac{1}{2\sqrt{x}}U_{2n-1}(\sqrt{x})\\
 & = & (2n-1)^{n-3}2^{(2n-1)(n-3)}\prod_{\substack{U'_{2n-1}(x_{i})=0\\
x_{i}>0
}
}\frac{4n^{2}-1}{2x_{i}(x_{i}^{2}-1)}U_{2n-1}(x_{i})\\
 & = & \frac{(4n^{2}-1)^{n-1}}{2^{n-1}\sqrt{2n/(2n-1)2^{2n-1}}U'_{2n-1}(1)/(2n-1)2^{2n-1}}\sqrt{\mbox{Disc}U_{2n-1}(x)}\\
 &  & \times\frac{(2n-1)^{n-3}2^{(2n-1)(n-3)}}{\sqrt{(2n-1)^{2n-1}2^{(2n-1)(2n-2)}}}\\
 & = & \frac{(4n^{2}-1)^{n-1}}{2^{n-1}\sqrt{2n/(2n-1)2^{2n-1}}U'_{2n-1}(1)/(2n-1)2^{2n-1}}\\
 &  & \times2^{(2n-1)^{2}/2}(2n)^{(2n-3)/2}\frac{\sqrt{2n-1}}{(2n-1)^{3}2^{2(2n-1)}}\\
 & = & (2n-1)^{n-2}(2n+1)^{n-1}n^{n-2}2^{n(2n-3)}/U'_{2n-1}(1),
\end{eqnarray*}
where $U'_{2n-1}(1)$ can be computed from its trigonometric definition:
\begin{eqnarray*}
U'_{2n-1}(1) & = & -\lim_{\theta\rightarrow0}\frac{2n\cos2n\theta\sin\theta-\sin2n\theta\cos\theta}{\sin^{3}\theta}\\
 & = & -\lim_{\theta\rightarrow0}\frac{2n(1-2n^{2}\theta^{2})(\theta-\theta^{3}/6)-(2n\theta-8n^{3}\theta^{3}/6)(1-\theta^{2}/2)}{\theta^{3}}\\
 & = & -\frac{2}{3}n(4n^{2}-1).
\end{eqnarray*}

Thus
\[
\mbox{Disc}_{x}U'_{2n-1}(\sqrt{x})=3(2n+1)^{n-2}(2n-1)^{n-3}n^{n-3}2^{2n^{2}-3n-1}.
\]
 Thus $\mbox{Disc}_{x}J(x)$ has factors 2, 3, and factors of powers
of $n$, $2n-1$, $2n+1$.

\section{The discriminant of $M(x)$}

In this section, we will present an explicit formula for the discriminant
of $M(x)$ and show that this discriminant also has factors 2, 3,
and factors of powers of $n$, $2n-1$, $2n+1$. In particular, we
will prove the theorem: 

\begin{theorem}

The discriminant of $M(x)$ is 
\[
\mbox{Disc}_{x}M(x)=\pm(2n-1)^{n-3}(2n+1)^{n-2}2^{2n^{2}-n-5}3n^{n-3}.
\]

\end{theorem}

We note that $M(x)$ is a polynomial in $x$ of degree $n-1$ whose
leading coefficient is $(2n-1)2^{2n}$. To compute $M'(x)$ we notice
the following fact:
\begin{equation}
(2n-1)T_{2n+1}(x)+(2n+1)T_{2n-1}(x)=2(2n-1)(2n+1)\int_{0}^{x}tU_{2n-1}(t)dt.\label{eq:sumT}
\end{equation}

From this, we have 
\begin{eqnarray*}
M'(x) & = & -\frac{3}{2}\frac{M(x)}{x}+\frac{(2n+1)(2n-1)}{x\sqrt{x}}U_{2n-1}(\sqrt{x}).
\end{eqnarray*}
Thus the definition of discriminant \eqref{eq:disdif} yields 
\begin{eqnarray*}
\mbox{Disc}_{x}(M(x)) & = & \pm(2n-1)^{n-3}2^{2n(n-3)}\\
 &  & \times(2n-1)^{n-1}(2n+1)^{n-1}\frac{(2n-1)2^{2n}}{M(0)}\prod_{M(x_{i})=0}x_{i}^{-1/2}U_{2n-1}(\sqrt{x_{i}}).
\end{eqnarray*}

where the value of the free coefficient of $M(x)$ can be obtained
from \eqref{eq:sumT}

\[
M(0)=\pm4(2n-1)(2n+1)n/3.
\]
 By substituting this value to the equation above, we obtain 
\[
\mbox{Disc}_{x}M(x)=\frac{\pm2^{2n^{2}-4n-2}(2n-1)^{2n-4}(2n+1)^{n-2}3}{n}\prod_{M(x_{i})=0}x_{i}^{-1/2}U_{2n-1}(\sqrt{x_{i}}).
\]

To compute the product above, we follow an idea from Jemal Gishe and
Mourad Ismail \cite{gi} by writing $M(x)$ as a linear combination
of $U_{2n-1}(x)$ and $U_{2n}(x)$. In particular this combination
can be obtained from \eqref{eq:TU1} and \eqref{eq:TU2} as below:
\begin{eqnarray*}
(2n-1)T_{2n+1}(\sqrt{x})+(2n+1)T_{2n-1}(\sqrt{x}) & = & -2T_{2n+1}(\sqrt{x})+2\sqrt{x}(2n+1)T_{2n}(\sqrt{x})\\
 & = & -U_{2n+1}(\sqrt{x})+U_{2n-1}(\sqrt{x})\\
 &  & +2\sqrt{x}(2n+1)(U_{2n}(\sqrt{x})-\sqrt{x}U_{2n-1}(\sqrt{x}))\\
 & = & 4n\sqrt{x}U_{2n}(\sqrt{x})-(2x(2n+1)-2)U_{2n-1}(\sqrt{x}).
\end{eqnarray*}
Thus
\begin{eqnarray*}
\prod_{M(x_{i})=0}x^{-1/2}U_{2n-1}(\sqrt{x}) & = & \frac{2^{(2n-1)(n-1)}}{(2n-1)^{n-1}2^{2n(n-1)}}\prod_{x_{i}^{-1/2}U_{2n-1}(\sqrt{x_{i}})=0}M(x_{i})\\
 & = & \frac{2^{(2n-1)(n-1)}}{(2n-1)^{n-1}2^{2n(n-1)}}4^{n-1}n^{n-1}\\
 &  & \times\frac{2^{2n-1}}{2n}\prod_{x^{-1/2}U_{2n-1}(\sqrt{x_{i}})=0}U_{2n}(\sqrt{x})\\
 & = & \frac{n^{n-2}2^{3n-3}}{(2n-1)^{n-1}}\prod_{\substack{U_{2n-1}(x_{i})=0\\
x_{i}>0
}
}U_{2n}(x_{i})\\
 & = & \pm\frac{n^{n-2}2^{3n-3}}{(2n-1)^{n-1}}
\end{eqnarray*}
where the last equality is obtained from \eqref{eq:resU}. Hence
\[
\mbox{Disc}_{x}M(x)=\pm(2n-1)^{n-3}(2n+1)^{n-2}2^{2n^{2}-n-5}3n^{n-3}.
\]
From this we conclude that the discriminants of $M(x)$ and $J(x)$
have the same factors 2, 3, and all the factors of $n$, $2n-1$,
$2n+1$.

\section{The roots of $M(x)$ and $J(x)$}

In this section we will show that the roots of $M(x)$ and $J(x)$
are pairwise similar except for one pair. In particular, we will prove
the theorem:

\begin{theorem}

Suppose $n$ is sufficiently large. Then for every root $x_{0}$ of
$J(x)$ except the smallest root, there is a root of $M(x)$ in the
interval
\[
[x_{0}-3/(10n^{2}),x_{0}+1/(2n^{2})).
\]

\end{theorem}

We first show that with one exception, the roots of $M(x)$ and $U'_{2n-1}(\sqrt{x})$
are pairwise close. From the definition, the roots of these two polynomial
are positive real numbers. To simplify the computations, we consider
$U'_{2n-1}(x)$ and the polynomial 
\[
R(x)=(2n-1)T_{2n+1}(x)+(2n+1)T_{2n-1}(x)
\]
whose roots (except 0) are the square roots of the positive roots
of $J(x)$ and $M(x)$ respectively. And it will suffice to consider
the positive roots of $R(x)$ and $U'_{2n-1}(x)$. 

Let $\zeta$ be a positive root of $U'_{2n-1}(x)$. We will show that
for a certain small value $\delta>0$, the quantities $R(\zeta)$
and $R(\zeta-\delta)$ have different signs and thus $R(x)$ admits
a root in the small interval $(\zeta-\delta,\zeta)$. First, we observe
that the equation \eqref{eq:rootJ} gives
\[
R(\zeta)=2U_{2n-1}(\zeta).
\]

For some $t\in(\zeta-\delta,\zeta)$ , we have

\[
R(\zeta-\delta)=R(\zeta)-R'(t)\delta.
\]
 where 
\[
R'(t)=2(2n+1)(2n-1)tU_{2n-1}(t).
\]

Thus
\begin{equation}
R(\zeta-\delta)=2U_{2n-1}(\zeta)-2(2n+1)(2n-1)tU_{2n-1}(t)\delta.\label{eq:taylor1}
\end{equation}

To prove $R(\zeta-\delta)$ and $R(\zeta)=2U_{2n-1}(\zeta)$ have
different signs, it remains to show that $2(2n+1)(2n-1)tU_{2n-1}(t)\delta$
is sufficiently large in magnitude and has the same sign as $R(\zeta)$.
To ensure this fact, we first impose the following two conditions
on $\delta$:

\begin{eqnarray}
\zeta\delta & > & \frac{A}{(2n+1)(2n-1)}\label{eq:con1}\\
\zeta & > & 6\delta\label{eq:con3}
\end{eqnarray}
where the value of $A$ and the existence of $\delta$ with respect
to these conditions will be determined later. With these conditions,
we obtain the following lower bound for $R'(t)\delta$: 
\begin{eqnarray}
|2(2n+1)(2n-1)U_{2n-1}(t)t\delta| & > & 2(2n+1)(2n-1)|U_{2n-1}(t)|(\zeta-\delta)\delta\nonumber \\
 & > & \frac{5}{3}(2n+1)(2n-1)|U_{2n-1}(t)|\zeta\delta\nonumber \\
 & > & \frac{5A}{3}|U_{2n-1}(t)|.\label{eq:lowerbound}
\end{eqnarray}
We now need to show that $U_{2n-1}(t)$ is not too small compared
to $U_{2n-1}(\zeta)$. The Taylor formula gives
\begin{equation}
U_{2n-1}(t)=U_{2n-1}(\zeta)+U''_{2n-1}(\zeta-\epsilon)(t-\zeta)^{2}\label{eq:taylor2}
\end{equation}
where $0<\epsilon<\zeta-t$. From \eqref{eq:difU}, we have, for $-1<x<1$,
the following trivial bound for $U'_{n}(x)$: 
\begin{equation}
|U'_{n}(x)|<\frac{2n}{1-x^{2}}.\label{eq:bounddifU}
\end{equation}
Thus the equation \eqref{eq:difj}, with $\sqrt{x}$ replaced by $x$,
implies that for large $n$ 
\begin{eqnarray}
|U''_{2n-1}(\zeta-\epsilon)| & \sim & \frac{4n^{2}}{1-(\zeta-\epsilon)^{2}}|U_{2n-1}(\zeta-\epsilon)|\label{eq:difdifU}\\
 & < & \frac{4n^{2}}{1-\zeta^{2}}|U_{2n-1}(\zeta)|.\nonumber 
\end{eqnarray}

We impose another condition on $\delta$: 
\begin{equation}
\delta<\frac{\sqrt{1-\zeta^{2}}}{4n}.\label{eq:con2}
\end{equation}

With this condition, we have 
\begin{eqnarray*}
|U''_{2n-1}(\zeta-\epsilon)(t-\zeta)^{2}| & < & \frac{1}{4}|U_{2n-1}(\zeta)|.
\end{eqnarray*}

Now the equation \eqref{eq:taylor2} implies that $U_{2n-1}(t)$ and
$U_{2n-1}(\zeta)$ have the same sign and 
\[
|U_{2n-1}(t)|>\frac{3}{4}|U_{2n-1}(\zeta)|.
\]
This inequality combined with \eqref{eq:taylor1} and \eqref{eq:lowerbound}
show that $R(\zeta-\delta)$ and $R(\zeta)$ have different signs
if 
\[
A>\frac{8}{5}.
\]
 Thus $R(x)$ has a root in the interval $(\zeta-\delta,\zeta).$ 

Next we will show that $\delta$ exists with respect to all the imposed
conditions and that its value is small as soon as $\zeta$ is not
the smallest positive root of $U'_{2n-1}(x)$. Let $\zeta$ be such
a root. From \eqref{eq:con1}, \eqref{eq:con3} and \eqref{eq:con2},
it suffices to show $\zeta$ satisfies the identity 
\begin{equation}
\frac{A}{(2n-1)(2n+1)}<\min\left(\frac{1}{4n}\zeta\sqrt{1-\zeta^{2}},\frac{\zeta^{2}}{6}\right).\label{eq:condidtion}
\end{equation}

The graph of the function $\zeta\sqrt{1-\zeta^{2}}$ is given below:

\begin{center}
\includegraphics[scale=0.5]{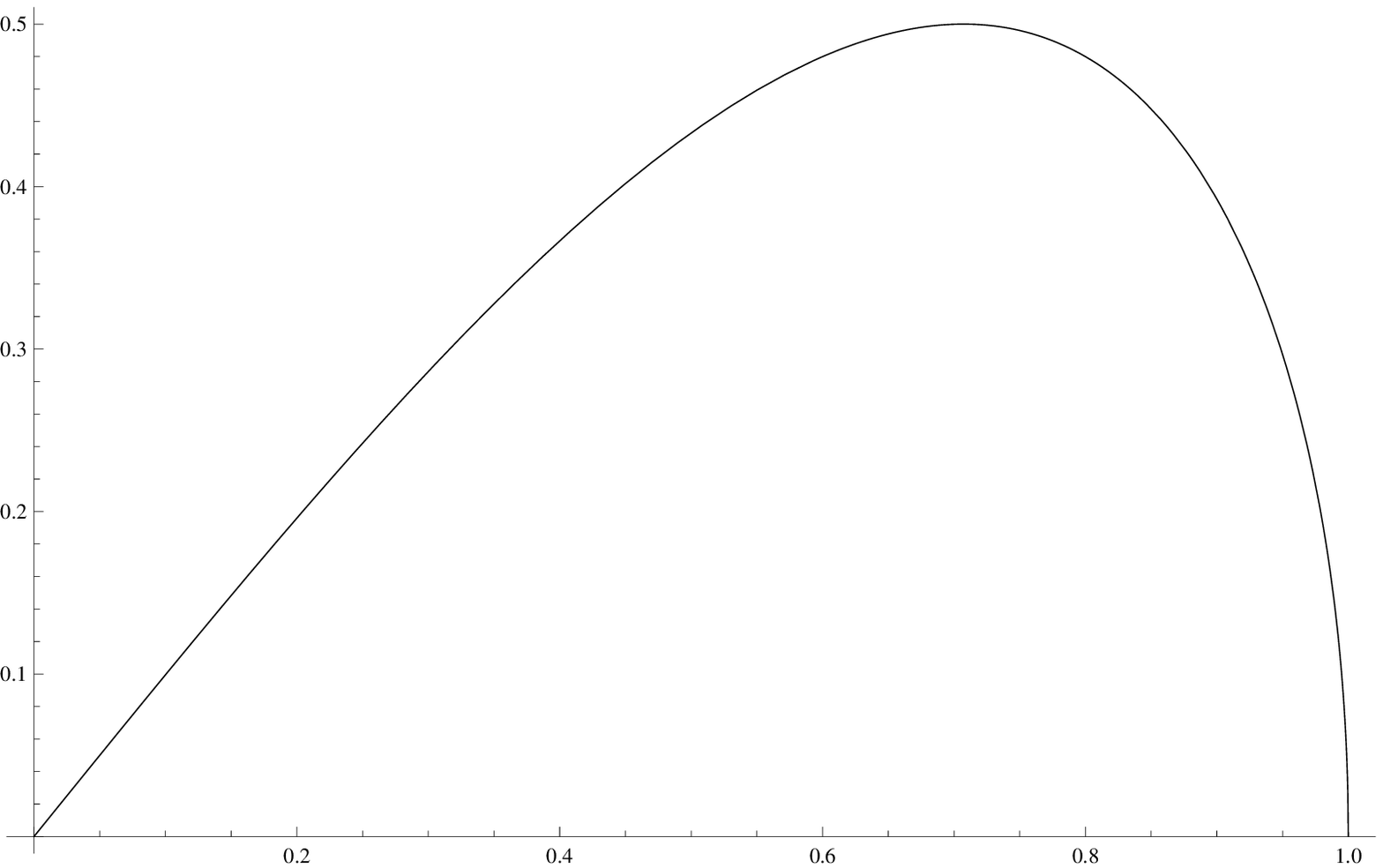}
\par\end{center}

From this graph, it suffices to show that $\zeta$ is not too small
and not too close to 1. Let $\zeta_{1}>0$ be the second smallest
root of $U'_{2n-1}(x)$. This root has an upper bound 
\[
\zeta_{1}<\cos\frac{\pi(n-2)}{2n}\sim\frac{\pi}{n}
\]
since the right side is the second smallest positive root of $U_{2n-1}(x)$.
The smallest positive root of $U_{2n-1}(x)$ is 
\[
\cos\frac{\pi(n-1)}{2n}\sim\frac{\pi}{2n}.
\]

Next the trivial lower bound $\zeta_{1}>\pi/2n$ does not guaranty
the existence of $A$ in \eqref{eq:condidtion}. We need a slightly
better lower bound. 

For some $\xi\in(\pi/2n,\zeta_{1})$ we have 
\[
(\zeta_{1}-\frac{\pi}{2n})U'_{2n-1}(\xi)=U_{2n-1}(\zeta_{1}).
\]
 The bound \eqref{eq:bounddifU} gives
\begin{eqnarray*}
\zeta_{1}-\frac{\pi}{2n} & > & \frac{|U_{2n-1}(\zeta_{1})|(1-\xi^{2})}{2(2n-1)}\\
 & > & \frac{1-(\pi/n)^{2}}{2(2n-1)}\\
 & \sim & \frac{1}{4n},
\end{eqnarray*}
so 
\[
\zeta_{1}>\frac{2\pi+1}{4n}.
\]

We now see that for $n$ large $\zeta_{1}$ satisfies the condition
\eqref{eq:condidtion} if 
\[
A<\frac{2\pi+1}{4}.
\]

Similarly, let $\zeta_{2}$ be the largest root of $U'_{2n-1}(x)$.
The trivial upper bound $\zeta_{2}<\cos\pi/2n$, the largest root
of $U_{2n-1}(x)$, does not guaranty the existence of $A$. However
for some $\xi\in(\zeta_{2},\cos\pi/2n)$ we have 
\[
(\cos\frac{\pi}{2n}-\zeta_{2})^{2}U''_{2n-1}(\xi)=U_{2n-1}(\zeta_{2}).
\]
The approximation \eqref{eq:difdifU} yields 
\begin{eqnarray*}
(\cos\frac{\pi}{2n}-\zeta_{2})^{2} & \sim & \frac{|U_{2n-1}(\zeta_{2})|(1-\xi^{2})}{4n^{2}|U_{2n-1}(\xi)|}\\
 & > & \frac{(1-\cos^{2}\pi/2n)}{4n^{2}}\\
 & \sim & \frac{\pi^{2}}{16n^{4}}.
\end{eqnarray*}
Thus 
\begin{eqnarray*}
\sqrt{1-\zeta_{2}^{2}} & > & \sqrt{1-(\cos\frac{\pi}{2n}-\frac{\pi}{4n^{2}})^{2}}\\
 & \sim & \sqrt{\frac{\pi^{2}}{4n^{2}}+\frac{\pi}{2n^{2}}}\\
 & > & \frac{2\pi+1}{4n}.
\end{eqnarray*}
Now the condition\eqref{eq:condidtion} is satisfied provided that
$A<(2\pi+1)/4$. 

Thus we can let $A$ be any value in the interval $(8/5,(2\pi+1)/4).$
Since for large $n$, $A$ and $\delta$ are arbitrarily bigger than
$8/5$ and $A/4n^{2}\zeta$ respectively, we have that for every root
$\zeta$ of $U'_{2n-1}(x)$ except the smallest root, there is a root
of $R(x)$ in the small interval 
\[
[\zeta-\frac{8}{20n^{2}\zeta},\zeta).
\]
This implies that in each pair, the root $\zeta^{2}$ of $U'_{2n-1}(\sqrt{x})$
is bigger than the root of $x_{i}$ of $M(x)$ by at most
\[
\frac{8}{20n^{2}\zeta}(\zeta+\sqrt{x_{i}})<\frac{8}{10n^{2}}.
\]
Combining this fact with \eqref{eq:rootJ}, we conclude that for every
root $x_{0}$ of $J(x)$ except the smallest root, there is a root
of $M(x)$ in the interval 
\[
[x_{0}-\frac{3}{10n^{2}},x_{0}+\frac{1}{2n^{2}})
\]
for large $n$. This concludes the proof.

Remark: Similar Mutt and Jeff phenomena seem to occur for the integral
transform of Legendre polynomials.

\end{document}